\documentclass{amsart}
\usepackage{amssymb}
\theoremstyle{plain}
\newtheorem{thm}{Theorem}[section]
\newtheorem{lem}[thm]{Lemma}
\newtheorem{que}[thm]{Question}
\newtheorem{cor}[thm]{Corollary}
\newtheorem{prop}[thm]{Propositon}
\theoremstyle{definition}
\newtheorem{defn}[thm]{Definition}

\newtheorem{remark}[thm]{Remark}

\begin{document}

\title[Detecting Codimension One Manifold Factors]{Detecting Codimension One Manifold Factors
with Topographical Techniques}

\author[D. M. Halverson and D. Repov\v s]{Denise M. Halverson and Du\v san Repov\v s}
\address{Department of Mathematics, Brigham Young University,
Provo, UT 84602}
\email{halverson@math.byu.edu}

\address{Faculty of Mathematics and Physics, and
Faculty of Education, University of Ljubljana,
P.O. Box 2964,
Ljubljana, Slovenia 1001}
\email{dusan.repovs@guest.arnes.si}

\date{\today}



\begin{abstract}
We prove recognition theorems for codimension one manifold factors
of dimension $n \geq 4$.  In particular, we formalize topographical
methods and introduce three ribbons properties: the crinkled ribbons
property, the twisted crinkled ribbons property, and the fuzzy
ribbons property. We show that $X \times \mathbb{R}$ is a manifold
in the cases when $X$ is a resolvable generalized manifold of finite
dimension $n \geq 3$ with either: (1) the crinkled ribbons property;
(2) the twisted crinkled ribbons property and the disjoint point
disk property; or (3) the fuzzy ribbons property.
\end{abstract}

\maketitle

\section{Introduction}

In this paper we provide general position techniques that fully
utilize a general position characterization of codimension one
manifold factors of dimension $n\geq 4$.
A {\it codimension one manifold
factor} is a space $X$ such that $X \times \mathbb{R}$ is a manifold.
The famous Cell-like Approximation Theorem of Edwards
\cite{Daverman book, Daverman-Halverson 2, Edwards 1, Edwards 2}
characterizes
the manifolds of dimension $n \geq 5$ as precisely the
finite-dimensional resolvable generalized manifolds with the disjoint disk
property.  In the same vein, it has been shown that codimension one
manifold factors of dimension $n \geq 4$ are precisely the
finite-dimensional resolvable generalized manifolds with the disjoint
concordances property.

However, up until now, practical methods of identifying spaces as
codimension one manifold factors have appealed to a weaker general
position property, the disjoint homotopies property. How to fully
utilize the disjoint concordances property has been somewhat
elusive. The ribbons properties introduced in this paper fulfill
this role. We  show that the ribbons properties, if satisfied by an
ANR $X$, will imply that $X$ has the disjoint concordances property
and hence $X \times \mathbb{R}$ has the disjoint disks property.
Therefore, a finite-dimensional resolvable generalized manifold $X$
is a codimension one manifold factor if it possesses one of the following: (1) the crinkled ribbons property; (2) the twisted
crinkled ribbons property and the disjoint point disk property;
or (3) the fuzzy ribbons property. For motivation see the surveys
\cite{HaRe}-\cite{Repovs3}.

\section{Manifold Factors and Characterizations}

As previously stated, a space $X$ is a \emph{codimension one
manifold factor} if $X \times \mathbb{R}$ is a manifold. The fact
that the codimension one manifold factors of finite dimension $n\geq
4$ are precisely the resolvable generalized manifolds $X$ such that
$X \times \mathbb{R}$ has the disjoint disks property follows as a
corollary of Edwards' Cell-like Approximation Theorem
(cf. \cite{Daverman book}).
Recall that a space $X$ is
said to be
\emph{resolvable} if there is a
manifold $M$ and a surjective
 map $f: M \to X$ which is cell-like (i.e., $f^{-1}(x)$ has the shape of a
point for all $x \in X$).  Moreover, $X$ is said to have the \emph{disjoint
disks property} (DDP) if every pair of maps $f,g: D^2 \to Y$ can be
approximated by maps that have disjoint images.

Edwards' Cell-like
Approximation Theorem states that the manifolds of dimension $n \geq
5$ are precisely the finite-dimensional resolvable generalized
manifolds with the disjoint disk property.  It is well known that
not all resolvable generalized manifolds of dimension $n \geq 5$
have the DDP (cf. \cite{Daverman book}). Thus not all resolvable generalized manifolds are
manifolds. (In dimension $\leq 2$ every generalized manifold is a topological manifold, whereas for the situation in dimensions 3 and 4 see \cite{Repovs1}-\cite{Repovs3}.)

In general, a space $X$ is said to satisfy the
\emph{$(m,n)$-disjoint disks property} ($(m,n)$-DDP) if any two maps
$f:D^m \to X$ and $g:D^n \to X$ can be approximated by maps with
disjoint images. As indicated previously, the $(2,2)$-DDP is simply
called the \emph{disjoint disks property} (DDP). The $(1,2)$-DDP is
called the \emph{disjoint arc-disk property} (DADP).  The
$(1,1)$-DDP is called the \emph{disjoint arcs property} (DAP). The
$(0,2)$-DDP is called the \emph{disjoint point-disk property} (DPDP).

All
generalized manifolds of dimension $n \geq 3$ are known to have the DAP.
A natural question is if the DADP, the middle dimension analogue of
the DAP and the DDP, provides a characterization of codimension one
manifold factors. As it turns out, the DADP condition is sufficient,
but not necessary, to determine if a finite-dimensional resolvable
generalized manifold of dimension $n \geq 4$ is a codimension one
manifold factor. Examples of codimension one manifold factors of
dimension $n \geq 4$ that fail to have the DADP can be found in \cite{Daverman
book, Daverman-Walsh,Halverson 3}.  In fact some of these examples
even fail to have the DPDP.

A list of general position properties that have proved useful in
recognizing codimension one manifold factors includes:

\begin{itemize}
  \item The disjoint arc-disk property \cite{Daverman 2}
  \item The disjoint homotopies property \cite{Halverson 3}
    \begin{itemize}
        \item The plentiful $2$-manifolds property \cite{Halverson 2}
        \item The method of $\delta$-fractured maps \cite{Halverson 2}
        \item The $0$-stitched disks property \cite{Halverson 3}
    \end{itemize}
  \item The disjoint concordances property \cite{Daverman-Halverson}
 \end{itemize}

It should be noted here that the disjoint concordances property is
the only property listed that provides a characterization of
codimension one manifold factors. Specifically, a resolvable
generalized manifold $X$ of finite dimension $n \geq 4$ is a
codimension one manifold factor if and only if $X$ satisfies the
disjoint concordances property.

\begin{defn} A \emph{path concordance} in a space $X$ is a
map $F:D \times I \to X \times I$ (where $D=I=[0,1]$) such that $F(D
\times e) \subset X \times e, e \in \{0,1\}.$  A metric space
($X,\rho$) satisfies the \emph{disjoint path concordances property
(DCP)} if, for any two path homotopies $f_i:D \times I \to X$
($i=1,2$) and any $\varepsilon > 0$, there exist path concordances
$F'_i: D \times I \to X \times I$ such that
\begin{center}
$F'_1(D \times I) \cap F'_2(D \times I) = \emptyset$
\end{center}
and $\rho (f_i, \text{proj}_X F'_i) < \varepsilon$.
\end{defn}
\noindent It is the main goal of this paper to establish practical
techniques that utilize this property.

In this paper we will be generalizing two properties: the plentiful
$2$-manifolds property and the method of $\delta$-fractures maps.
These properties were developed specifically to detect the disjoint
homotopies property in certain settings.  We will demonstrate how
the analogous ribbons properties can be used to detect the weaker
disjoint concordances property.

\section{Topographies}

We begin by restating the disjoint concordances property from a more
functional perspective.  In particular, we will restate the disjoint
concordance property in terms of the topographies.

\begin{defn}
A \emph{topography $\Upsilon$ on $Z$} is a partition of $Z$ induced
by a map $\tau: Z \to I$.  The \emph{$t$-level of $\Upsilon$} is
given by $$ \Upsilon_t = \tau^{-1}(t).$$
\end{defn}

\begin{defn}
A \emph{topographical map pair} is an ordered pair of maps $(f,
\tau)$ such that $f:Z \to X$ and $\tau: Z \to I$.  The map $f$ will
be referred to as the \emph{spatial map}
and
the map $\tau$ will be
referred to as the \emph{level map}. The topography associated with
$(f, \tau)$ is $\Upsilon$,
where $ \Upsilon_t = \tau^{-1}(t).$
\end{defn}

Note that a homotopy $f: Z\times I \to X$ has a naturally associated
topography, where $\tau: Z\times I \to I$ is defined by $\tau(x,t) =
t$. In particular, we may view $f: Z\times I \to X$ as being equivalent to
$(f,\tau)$ and we will refer to $(f,\tau)$ as the \emph{natural
topographical map pair associated with $f$}.

\begin{defn}
Suppose that for $i=1,2$, $\Upsilon^i$ is a topography on $Z_i$
induced by $\tau_i$ and $f_i: Z_i \to X$.  Then $(f_1, \tau_1)$ and
$(f_2, \tau_2)$ are \emph{disjoint topographical map pairs} provided
that for all $t \in I$,
$$f_1( \Upsilon^1_t) \cap f_2(\Upsilon^2_t) = \emptyset.$$
A space $X$ has the \emph{disjoint topographies property} if any two
topographical map pairs $(f_i, \tau_i)$ ($i=1,2$), where $f_i: D^2
\to X$, can be approximated by disjoint topographical map pairs.
\end{defn}

The proof of the
following result is straightforward:

\begin{thm}
An ANR $X$ has the disjoint topographies property if and only if $X
\times \mathbb{R}$ has the disjoint disks property.
\end{thm}

\begin{proof}
Suppose $X$ has the disjoint topographies property.  For $i=1,2$,
let $F_i: D^2 \to X \times I$.  Let $\text{proj}_X: D^2 \to X$ and
$\text{proj}_I: D^2 \to I$ be the standard projection maps. Define
$f_i = \text{proj}_X \circ F_i$ and $\tau_i = \text{proj}_I \circ
F_i$. Applying the disjoint topographies property we get disjoint
topographical map pairs $(f'_i, \tau'_i)$ that are approximations of
$(f_i, \tau_i)$. Then $F'_i = f'_i \times \tau'_i$ are the desired
approximations of $F_i$ with disjoint images.

Suppose that $X \times \mathbb{R}$ has the disjoint disks property.
Let $(f_i, \tau_i)$ be topographical map pairs for $i=1,2$.  Then
$F_i = f_i \times \tau_i: D^2 \to X \times I$.  By the disjoint
disks property $F_i$ can be approximated by $F'_i$ with disjoint
images. Let $f'_i = \text{proj}_X \circ F'_i$ and $\tau'_i =
\text{proj}_I \circ F'_i$. Then $(f'_i, \tau'_i)$ are the desired
disjoint topographical map pairs approximating $(f_i, \tau_i)$.
\end{proof}

\noindent However, this result is not the main focus of this paper.
Our aim is to provide alternative equivalent conditions which are
more easily verified.

\section{Special Category Approximation Properties}

We desire to more carefully investigate the disjoint topographies
property so as to give it practical utility.  Similar to the
disjoint homotopies property analyzed in \cite{Halverson 1},
 the question of whether a space has the disjoint
concordances property ultimately reduces to the following question:
given a constant homotopy of a $1$-complex and an arbitrary homotopy
on another $1$-complex, can the natural topographical map pairs
associated with these homotopies be adjusted with ``control'' so as
to form disjoint topographical map pairs?   In this section, we will
clarify these characterizing conditions. In Section \ref{Equiv Thm}
we will demonstrate that the conditions give the desired result. The
ribbons properties in Sections \ref{Crinkled} and \ref{Fuzzy} will
specify practical circumstances in which these conditions may be
obtained.

\begin{defn}
A topographical map pair $(f, \tau)$ is in the \emph{$\mathcal{Z}$
category} if $f:Z \times I \to X$ and $\tau: Z \times I \to I$ so
that $Z \times \{e\} \subset \tau^{-1}(e)$ for $e=0,1$. We denote
$(f, \tau) \in \mathcal{Z}$.
\end{defn}

The \emph{$\mathcal{D}$ category} is defined by letting $Z = D =
[0,1]$. The \emph{$\mathcal{K}$ category} is defined by letting $Z =
K$, for some $1$-complex.

\begin{defn}
A topographical map pair $(f,\tau)$ is in the \emph{$\mathcal{Z}_c$
category} if

\begin{enumerate}
\item $(f,\tau) \in \mathcal{Z}$;
\item $f:Z \times I \to X$ is a constant
homotopy;
and
\item $(f,\tau)$ is the natural topographical map pair associated with $f$.
\end{enumerate}
\end{defn}

For emphasis on the relevant characteristics, we define the
conditions that will be the main focus of the next section in two
stages.  In the following definitions, the notation $\mathcal{Z}_i$
is intended to represent a category such as $\mathcal{D}$ or
$\mathcal{K}$.

\begin{defn} A space $X$ has the \emph{$\mathcal{Z}_1 \times
\mathcal{Z}_2$ category disjoint topographies property (
$\mathcal{Z}_1 \times \mathcal{Z}_2$ DTP)} if any two topographical
map pairs $(f_i, \tau_i) \in \mathcal{Z}_i$, for $i = 1,2$, can be
approximated by disjoint topographical map pairs $(f'_i, \tau'_i)
\in \mathcal{Z}_i$.
\end{defn}

\begin{defn} A space $X$ has the \emph{$\mathcal{Z}_1 \times
\mathcal{Z}_2$  DTP*} if for any pair of maps $(f_i, \tau_i) \in
\mathcal{Z}_i$, for $i = 1,2$, there are maps $(f'_i, \tau'_i) \in
\mathcal{Z}_i$ so that each $f_i'$ is an approximation of $f_i$.
\end{defn}

\noindent Note specifically that the
\emph{$\mathcal{Z}_1 \times \mathcal{Z}_2$  DTP*} condition does not require
the maps $\tau'_i$ to approximate $\tau_i$.

A careful look at the definitions will reveal that the $\mathcal{D}
\times \mathcal{D}$ DTP* is just the disjoint concordance property
in the language of topographies. Our goal will be to show that the
disjoint concordances property is equivalent to more versatile
conditions, namely the $\mathcal{K}_c \times \mathcal{K}$ DTP* and
the $\mathcal{D}_c \times \mathcal{D}$ DTP* in the case that the target
space of the spatial map has the $(0,2)$-DDP. It is these conditions to
which our ribbons properties appeal.

\section{Extension Theorems}

In this section we recall a couple of classical extension theorems
that are used extensively when performing general position
adjustments in ANR's.  We also establish specific extension theorems
applicable in the setting of spaces with the various disjoint
topographies properties. To see a proof of the following homotopy
extension theorem the reader can refer to \cite{Halverson 1}.

\begin{thm} [Homotopy Extension Theorem (HET)]
Suppose that $f:Y \to X$ is a continuous map where $Y$ is a metric
space and $X$ is an ANR, $Z$ is a compact subset of $Y$ and
$\varepsilon > 0$. Then there exists $\delta > 0$ such that each
$g_Z:Z \to X$ which is $\delta$-close to  $f|_Z$ extends to $g:Y \to
X$ so that $g$ is $\varepsilon$-homotopic to $f$. In particular, for
any open set $U$ such that $Z \subset U \subset Y$, there is a
homotopy $H:Y \times I \to X$ so that:
\begin{enumerate}
    \item $H_0 = f$ and $H_1 = g$;
    \item $g|_Z = g_Z$;
    \item $H_t|_{Y-U} = f|_{Y-U}$, for all $t \in I$; and
    \item $diam(H(y \times I)) < \varepsilon$ for all $y \in Y$.
\end{enumerate}
\label{HET}
\end{thm}

\begin{cor}[Map Extension Theorem (MET)] \label{MET}
Suppose that $f:Y \to X$ is a continuous map where $Y$ is a metric
space and $X$ is an ANR, $Z$ is a compact subset of
$Y$ and
$\varepsilon > 0$.  Then there exists $\delta > 0$ such that each
$g_Z: Z \to X$ which is $\delta$-close to $f|_Z$ extends to $g:Y \to
X$ so that $\rho(f,g)< \varepsilon$.
\end{cor}

\noindent In the arguments that follow, when we say that ``without
loss of generality (such and such) maps into an ANR are already
adjusted to exploit (some) general position property'', we are
generally appealing to an application of MET.  For example, given
maps $f_i: D^2 \to X$, i=1,2, where $X$ is an ANR, with the DAP when we
say that we may assume without loss of generality that the restrictions
of these maps to a finite (or countable) collection of arcs in the
domain have disjoint images, we are applying MET.

\begin{cor}[Special DTP Extension Theorem]
Let $X$ be an ANR. Suppose that for $i=1,2$, $(f_i, \tau_i)$ are
topographical map pairs so that $f_i: Y_i  \to X$ and $\tau_i: Y_i
 \to I$, where $Y_i$ is a compact metric space. Suppose
further that $A_i \subset Y_i $ is compact so that:
\begin{enumerate}
\item $(f_1 |_{A_1}, \tau_1 |_{A_1})$ and $(f_2, \tau_2)$ are
disjoint topographical map pairs;
\item $(f_1, \tau_1)$ and $(f_2 |_{A_2}, \tau_2 |_{A_2})$ are
disjoint topographical map pairs; and
\item $(f_1, \tau_1)$ and $(f_2, \tau_2)$ can be approximated by disjoint
topographical map pairs.
\end{enumerate}
Then $(f_i,\tau_i)$ can be approximated by disjoint topographical
map pairs $(f_i',\tau_i')$ so that $(f_i' |_{A_i}, \tau_i' |_{A_i})
= (f_i |_{A_i}, \tau_i |_{A_i})$.
\end{cor}

\begin{proof}
Suppose the objects in the hypothesis are given.  By continuity and
local compactness we can find compact neighborhoods
$N_i$ of
$A_i$ so that (1) and (2) still hold when $A_1$ is replaced with
$N_i$.  Choose $\varepsilon
> 0$ so that (1) and
(2) still hold with $N_i$ replaced with $A_i$  and $(f_i,\tau_i)$
replaced with any $\varepsilon$-approximation of $(f_i,\tau_i)$.

Let
$\delta_i>0$ be values promised by the HET for $(f_i, \tau_i)$ and
choose $\delta > 0$ so that $\delta < \delta_1, \delta_2$. Find
$\delta$-approximations $(f'_i,\tau'_i)$ of $(f_i,\tau_i)$ that are
disjoint topographical map pairs. Let $Z_i = \overline{Y_i - N_i}$,
$U_i = Y_i-A_i$, and $g_i = f'_i |_{Z_i}$. Let $f'_i:Y_i \to X$ be
the end of the homotopy $H^i: Y_i \times I \to X$ promised by the HET.
Then $(f'_i, \tau'_i)$ are the desired disjoint homotopies such that
$(f_i' |_{A_i}, \tau_i' |_{A_i}) = (f_i |_{A_i}, \tau_i |_{A_i})$.
\end{proof}

In the next result, the \emph{end levels of a concordance $F: Y
\times I \to X \times I$ or  a topographical map pair $(f,\tau)$}
defined on $Y \times I$ will refer to $Y \times \{0 \}$ and $Y\times
\{1\}$.

\begin{prop} \label{ext2}
Let $X$ be an ANR.  Suppose $(f_i, \tau_i) \in \mathcal{Z}$ such
that $f_i: Y_i \times I \to X$.  If the restriction to the end
levels is a disjoint topographical map pair and $(f_i, \tau_i)$ can
be approximated by disjoint topographical map pairs, then $(f_i,
\tau_i)$ can be approximated by topographical map pairs fixed on the
end levels. An analogous result is true for concordances.
\end{prop}

\begin{proof}
Let $E_i = Y_i \times \{0,1\}$. By hypothesis, $(f_i|_{E_i}, \tau_i
|_{E_i})$ are disjoint topographical map pairs. Let $\varepsilon
> 0$ so that any $\varepsilon$-approximation of $(f_i|_{E_i}, \tau_i
|_{E_i})$ are still disjoint topographical map pairs.  Let
$\delta_i> 0$ be a value promised by the MET for $\varepsilon$ and
$f_i|_{E_i}$.  Choose $\delta> 0$ so that $\delta < \delta_1,
\delta_2$. Let $(g_i, \mu_i)$ be $\delta$-approximations of $(f_i,
\tau_i)$ in $\mathcal{Z}$. Then there are $\varepsilon$ homotopies
between $(f_i|_{E_i}, \tau_i |_{E_i})$ and $(g_i|_{E_i}, \mu_i
|_{E_i})$, call these $H^i: E_i \times I \to X$.

For $0<\zeta<\frac{1}2$, let $\theta_{\zeta}: [0,1] \to [\zeta,
1-\zeta]$ be the standard order preserving linear  map.  Define
$f_i^{\zeta}: Y_i \times I \to X$ such that:

$$f_i^{\zeta}(x,t) = \left\{
    \begin{array}{ll}
        H^i((x,0),\frac{t}{\zeta}) & \hbox{ if } t \in [0, \zeta) \\
      g_i(x,\theta^{-1}_{\zeta} (t) ) \, & \hbox{ if } t \in [ \zeta, 1-\zeta]\\
      H^i((x,1),\frac{1-t}{\zeta}) & \hbox{ if } t \in ( 1-\zeta, 1]
    \end{array}
  \right.
$$
and define
$\tau_i^{\zeta}: Y_i \times I \to X$ such that:
$$\tau_i^{\zeta}(x,t) = \left\{
    \begin{array}{ll}
     \theta _{\zeta}\mu_i(x,t) \, & \hbox{ if } t \in [ \zeta, 1-\zeta]\\
      t & \hbox{ if } t \in [0, \zeta) \cup ( 1-\zeta, 1]
    \end{array}
  \right.
$$

\noindent Note that $(f^{\zeta}_i |_{E_i},
\tau^{\zeta}_i|_{E_i})=(f_i|_{E_i}, \tau_i|_{E_i})$. Moreover, for
sufficiently small $\zeta$, $(f^{\zeta}_i, \tau^{\zeta}_i)$ are
disjoint topographical map pairs that are an $\varepsilon$-approximation of $(f_i, \tau_i)$.

The argument for concordances is analogous.
\end{proof}

\section{Tools for Finding Disjoint Topographies}

The following four ``R'' strategies can be used to manipulate
topographical map pairs to be disjoint:

\begin{enumerate}
\item {\it Reimage} - modify the spatial image set;
\item {\it Realign} - modify the position of the levels by a self homeomorphism of the domain $D
\times I$;
\item {\it Reparametrize} - relabel the levels by a continuous map fixing
the $t=0$ and $t=1$ levels; and
\item {\it  moRph} - redefine the topographical structure.
\end{enumerate}

The first strategy is realized by adjusting the spatial maps.  The
last three are realized by adjusting the level maps.
It is the fourth strategy that is unique to the topographies
approach, adding flexibility in that the shape of the levels can be
changed.  This is the maneuver that puts the topographical approach
at an advantage over the homotopies approach in detecting
codimension one manifold factors.  It is this last strategy that
will be fully exploited by the new ribbons properties of Sections
\ref{Crinkled} and \ref{Fuzzy}.

This section will be devoted to adapting several basic tools that
are useful in constructing approximating disjoint topographical map
pairs. The first two results are generalizations of results obtained
for homotopies found in \cite{Halverson 1}.

\begin{defn}
Suppose for $i=1,2$ that $(f_i, \tau_i)$ are topographical map pairs
having topographies $\Upsilon^i$.  Then the \emph{set of
parameterization points of intersection}, denoted by
PPIN$((f_1,
\tau_1),(f_2, \tau_2))$ is
$$\text{PPIN}((f_1, \tau_1),(f_2, \tau_2)) = \{ (t_1, t_2) \in I^2 \ | \ f_1(\Upsilon^1_{t_1}) \cap f_2(\Upsilon^2_{t_2}) \ne
\emptyset \}.$$
\end{defn}

In the next result we show that if PPIN$((f_1, \tau_1),(f_2,
\tau_2))$  is $0$-dimensional, then we may obtain approximating
disjoint topographical map pairs by reparametrizing the levels of
the topography.  In particular, a reparametrization is a relabeling
of the levels determined by replacing $t$ with a function
$\gamma(t)$.

\begin{lem}[Reparametrization Lemma]  Suppose
for $i=1,2$ that
$(f_i, \tau_i)$ are topographical map pairs having topographies
$\Upsilon^i$ such that PPIN$((f_1, \tau_1),(f_2, \tau_2))$ is
$0$-dimensional and $f_1(\Upsilon^1_e) \cap f_2(\Upsilon^2_e) =
\emptyset$, for $e=0,1$. Then there are arbitrarily close
approximations $\tau'_i$ of $\tau_i$ so that $(f_1, \tau'_1)$ and
$(f_2, \tau_2')$ are disjoint topographical map pairs.
\end{lem}

\begin{proof}
Suppose $\varepsilon
> 0$. Since $Z= \text{PPIN}((f_1, \tau_1),(f_2, \tau_2))$ is
$0$-dimensional there is a path $\gamma: I \to I \times I - Z$ from
$(0,0)$ to $(1,1)$ such
that $| \gamma(t) - t | < \varepsilon$. Let
$\tau'_i = \gamma \circ \tau_i$. Then $(f_i, \tau'_i)$ are disjoint
topographical map pairs.
\end{proof}

\begin{prop}  Suppose $X$ is a locally compact ANR with the DAP.  Then $X$ has
the disjoint topographies property if and only if for any pair of
topographical maps $(f_i, \tau_i)$, for i=1,2 such that $f_i: D
\times I \to X$, there exist arbitrarily close approximations
$(f'_i, \tau'_i)$, such that PPIN$((f'_1, \tau'_1),(f'_2, \tau'_2))$
is $0$-dimensional.
\end{prop}

\begin{proof}
To see the forward direction, assume
without loss of generality,
that $$f_1(D \times \{0,1\} \cup \{0,1\} \times I) \cap f_2(D \times
\{0,1\} \cup \{0,1\} \times I) = \emptyset$$   where $D \times I =
[0,1] \times [0,1]$.  It follows from the hypothesis that the
collection of maps $(f'_1, \tau'_1,f'_2, \tau'_2)$ such that $(f_1',
\tau_1')$ and $(f_2', \tau_2')$ are disjoint topographical map pairs
is dense in $\mathcal{D} \times \mathcal{D}$ and this collection is
clearly open by continuity arguments. Let $\gamma_k: I \to I$ be a
countable collection of maps such that the complement of the images
in the interior of $I \times I$ is $0$-dimensional and $\gamma_k(e)
= e$ for $e=0,1$. Find approximations $(f'_i, \tau'_i)$ so that
$(f'_1, \tau'_1)$ and $(f'_2, \gamma_k \tau'_2)$ are disjoint
topographical map pairs for all $k$.  Then PPIN$((f'_1,
\tau'_1),(f'_2, \tau'_2))$ is $0$-dimensional.

The reverse direction follows almost immediately from the
Reparametrization Lemma. The only technicality is that we need to
satisfy $\Upsilon^1_e \cap \Upsilon^2_e = \emptyset$. We may modify
$(f_i, \tau_i)$ by assuming that $\tau_i$ is a piecewise linear general position
map with care taken so that $\Upsilon^i_e = D \times \{e\}$ for $e
=0,1$.  Then we apply the DAP and the MET to modify $f_i$ so that
$f_1(\Upsilon^1) \cap f_2(\Upsilon^2) = \emptyset$.
\end{proof}

\begin{prop}\label{point arc homotopy}  Let $X$ be a finite-dimensional
ANR with the $(m_1 - 1, m_2)$-DDP and the $(m_1, m_2-1)$-DDP.
Suppose that $(f_i, \tau_i)$ are topographical map pairs such that
$f_i: Y_i \to X$ and $\tau_i: Y_i \to X$, where  $Y_i$ is a
$k$-complex such that $k \leq m_i$. Then there exist
approximations
$(f'_i, \tau'_i)$ that are disjoint topographical map pairs.
\end{prop}

\begin{proof}
Begin by modifying the maps $\tau_i: Y_i \to I$, if necessary, so
that each level is a $(k-1)$-complex. This can be accomplished by
approximating $\tau_i$ by a piecewise linear map in general
position. Next, apply the $(m_1 - 1, m_2)$-DDP and the $(m_1,
m_2-1)$-DDP conditions to adjust the maps $f_i$ so that each
rational level of $f_1$ is disjoint from the image of $f_2$ and each
rational level of $f_2$ is disjoint from the image of $f_1$. Denote
the adjusted maps by $(f'_i, \tau'_i)$. Then PPIN$((f'_1,
\tau'_1),(f'_2, \tau'_2))$ is closed $0$-dimensional. If follows by
the Reparametrization Lemma that $(f'_i, \tau'_i)$ can be
approximated by disjoint topographical map pairs.
\end{proof}

\section{Equivalence Theorem} \label{Equiv Thm}

In this section, we will demonstrate the following equivalence
theorem:

\begin{thm}[Equivalence Theorem]
Let $X$ be a locally compact separable ANR with the DAP.  Consider the
statements:
\begin{enumerate}
\item[(a)] $X$ has the $\mathcal{D}_c \times \mathcal{D}$ DTP*.
\item[(b)] $X$ has the $\mathcal{K}_c \times \mathcal{K}$ DTP*.
\item[(c)] $X$ has the $\mathcal{D} \times \mathcal{D}$ DTP*.
\item[(d)] $X$ has the disjoint concordance property.
\item[(e)] $X \times \mathbb{R}$ has the disjoint disks property.
\end{enumerate}
Then (b)-(e) are equivalent.  If in addition, $X$ has the $(0,2)$-DDP,
then (a)-(e) are equivalent.
\end{thm}

\begin{proof}
 Observe that (c) and (d) are trivially equivalent since the $\mathcal{D}
\times \mathcal{D}$ DTP* is the disjoint concordance property in the
language of topographies. In particular, equate $(f_i, \tau_i)$ as a
topographical map pair with $F_i = f_i \times \tau_i$ as a
concordance.  The fact that (d) and (e) are equivalent was the main
result established in \cite{Daverman-Halverson 2}. The fact that (c)
implies (a) is trivial since $\mathcal{D}_c \subset \mathcal{D}$.

It suffices to show that: (e) implies (b); (b) implies (c); and (a)
implies (b) in the case that $X$ has the $(0,2)$-DDP.

\medskip{ \bf (e) $\mathbf{\implies}$ (b)}:  In a locally compact separable ANR, the DDP condition is equivalent
to having the property that any two maps $\lambda_i: P_i \to X$ can
be approximated by maps with disjoint images where $P_i$ are
$2$-complexes (see \cite[Proposition 24.1]{Daverman book}). Given
topographical map pairs $(f_1, \tau_1) \in \mathcal{K}_c$ and $(f_2,
\tau_2) \in \mathcal{K}$ we may assume without loss of generality,
by applying the DAP, that the restrictions to the end levels are disjoint
topographical map pairs.  By Proposition \ref{ext2}, there are
approximations $(f'_i, \tau'_i)$ that are disjoint topographical map
pairs fixed on the end levels.  These are the desired
approximations.

\medskip{ \bf (b) $\mathbf{\implies}$ (c)}:  Let $(f_i, \tau_i) \in \mathcal{D}$ for $i=1,2$.
According to \cite[Theorem 3.3]{Halverson 1} there are piecewise
linear approximations  $\tau'_i$ of $\tau_i$ such that there exist:
\begin{enumerate}
\item a collection of $1$-complexes $K^i_1, \ldots, K^i_n$; and
\item a collection maps $\phi^i_j: K^i_i \times [t^i_{j-1},t^i_j] \to
D^2$ so that:
    \begin{enumerate}
    \item $\bigcup \text{im} (\phi^i_j )= D^2$;
    \item $\phi^i_j: K^i_j \times [t^i_{j-1},t^i_j] \to D^2$ is an embedding away from $K^i_j \times \{t^i_{j-1},t^i_j\}$; and
    \item $\tau'_i \circ \phi^i_j$ is a level preserving map.
    \end{enumerate}
\end{enumerate}
Without loss of generality we may assume that $t^1_j = t^2_j$ by
subdividing into smaller intervals if necessary. Thus we will denote
$t_j = t^1_j = t^2_j$.

Denote $L^i_j = \tau_i^{-1}(t_j)$.  These $1$-complexes are called
the transition levels.  Note that $$L_0^i = \phi^i_{1} (K^i_{1}
\times \{ t_0 \}) \text{ and } L_n^i = \phi^i_n (K^i_n \times \{ t_n
\})$$ and for $j=1, \ldots, n-1$,
$$L^i_j=\phi^i_j (K^i_j \times \{ t_j \}) \cup \phi^i_{j+1}
(K^i_{j+1} \times \{ t_j \}).$$   By applying the DAP, we may assume
that $f_1 (L^1_j) \cap f_2 (L^2_j) = \emptyset$.

Without loss of generality we may assume that the subintervals
$[t_{j-1}, t_j]$ are sufficiently small so that the adjustments that
will now follow will also be small. We will begin by modifying
the maps $\phi^i_j$ so that the $\mathcal{K}_c \times \mathcal{K}$
DTP* condition may be exploited. For $j = 1, \ldots n$, let $s_j =
\dfrac{t_{j-1} + t_j}{2}$. Define maps:

\begin{align}
\theta^1_j:& K^1_j \times [t_{j-1}, s_j] \to D^2; \theta^1_j(z,t) =
\phi^1_j(z, 2t-t_{j-1}) \\
\lambda^1_j:& K^1_j \times [s_j, t_j] \to D^2; \lambda^1_j(z,t) =
\phi^1_j(z,t_j) \\
\theta^2_j:& K^2_j \times [t_{j-1}, s_j] \to D^2; \theta^2_j(z,t) =
\phi^2_j(z,t_{j-1}) \\
\lambda^2_j:& K^2_j \times [s_j, t_j] \to D^2; \lambda^2_j(z,t) =
\phi^2_j(z, 2t-t_j)
\end{align}

Consider the natural topographical map pairs $(f_1\theta^1_j,
\eta^1_j)$, $(f_1\lambda^1_j, \mu^1_j)$, $(f_2\theta^2_j,
\eta^2_j)$, and $(f_2\lambda^2_j, \mu^2_j)$ defined for these
homotopies, respectively.  By applying the $\mathcal{K} \times
\mathcal{K}_c$ DTP* condition and Proposition \ref{ext2} on
$(f_1\theta^1_j, \eta^1_j)$ and $(f_2\theta^2_j, \eta^2_j)$ we can
find approximations $(g^1_j, \tilde{\eta}^1_j)$ and $(g^2_j,
\tilde{\eta}^2_j)$, respectively, that are disjoint topographical
map pairs fixed on the end levels. Likewise, there are
approximations of $(f_1\lambda^1_j, \mu^1_j)$ and $(f_2\lambda^2_j,
\mu^2_j)$, namely $(h^1_j, \tilde{\mu}^1_j)$ and $(h^2_j,
\tilde{\mu}^2_j)$, respectively, that are disjoint topographical map
pairs also fixed on the end levels.

Let
$$ f'_i(x) = \left\{
     \begin{array}{ll}
       g^i_j(\phi^i_j)^{-1}(x) & \hbox{ if } \tau'(x) \in [ t_{i-1}, s_i]\\
       h^i_j(\phi^i_j)^{-1}(x) & \hbox{ if } \tau'(x) \in [s_i, t_i]
     \end{array}
   \right.$$
and
$$ \tau'_i(x) = \left\{
     \begin{array}{ll}
       \tilde{\eta}^i_j(\phi^i_j)^{-1}(x) & \hbox{ if } \tau'(x) \in [ t_{i-1}, s_i]\\
       \tilde{\mu}^i_j(\phi^i_j)^{-1}(x) & \hbox{ if } \tau'(x) \in [s_i,
t_i].
     \end{array}
   \right.$$

\noindent Then $(f'_i, \tau'_i)$ are the desired disjoint
topographical map pairs in $\mathcal{D}$ that are approximations of
$(f_i, \tau_i)$.

\medskip{ \bf (a) $\mathbf{\implies}$ (b)}:  Note that this is the only case that requires the $(0,2)$-DDP condition.  Let $(f_1, \tau_1) \in \mathcal{K}_c$ and
$(f_2, \tau_2) \in \mathcal{K}$, where $f_i: K_i \times I \to X$ and
$\tau_i: K_i \times I \to I$.  Let $A_i$ be the $1$-complex $[K_i
\times \{0,1\}] \cup [K_i^{(0)} \times I]$. By applying the DAP and
Corollary \ref{point arc homotopy} (using the DAP and the $(0,2)$-DDP) we
may assume without loss of generality that:

\begin{enumerate}
\item $(f_1 | A_1, \tau_1|A_1)$ and $(f_2,\tau_2)$ are disjoint
topographical map pairs;
\item  $(f_1,\tau_1)$ and $(f_2 | A_2, \tau_2|A_2)$ are disjoint
topographical map pairs;
\item the restriction of $f_1$ to $K_1 \times \{0\}$ is an embedding;
\item the restriction $f_2$ to $A_2$ is an embedding; and
\item $f_1(A_1) \cap f_2(A_2) = \emptyset$.

\end{enumerate}

We wish to define topographical maps $(g_i, \eta_i) \in \mathcal{D}$
that can guide the appropriate modifications of $(f_i, \tau_i)$ to
give approximations that are disjoint topographical map pairs.  To
this end, using also the DAP, we may find $\alpha_i: D \to X$ so that:

\begin{enumerate}
\item[(6)] $f_i(K_i \times \{0\}) \subset \alpha_i(D)$;
\item[(7)] $\alpha_1(D) \cap \alpha_2(D) = \emptyset$; and
\item[(8)] $\alpha_i$ is a piecewise linear embedding, and in particular the restriction of $\alpha_i$
to $\alpha_i^{-1}f_i\left((\sigma- \sigma^{(0)}) \times \{ 0
\}\right) $ is an embedding for each simplex $\sigma \in K_i$.
\end{enumerate}

The maps $\alpha_i$ will determine the $0$-level maps of the new
maps $g_i: D \times I \to X$ which we will now construct.  For
reference in $D$, let $Q_i = \alpha_i^{-1}(f_i(K_i \times \{0 \}))$
and $P_i = \alpha_i^{-1}(f_i(K_i^{(0)} \times \{0 \}))$. For
reference in $D \times I$, define $B_i = Q_i \times\{0,1\} \cup P_i
\times I$ and $C_i = Q_i \times I$.

Define $g_1: D \times I \to X$
to be the constant homotopy so that $g_1(x,t) = \alpha_1(x)$. Let
$\eta_1: D \times I \to I$ be the standard projection map. Define
$\tilde{g}_2: D \times \{0\} \cup C_2 \to X$ so that:
$$\tilde{g}_2(x,t) = \begin{cases}
\alpha_2(x) & \text{ if } t=0 \\
f_2(z,t) & \text{ if } x \in Q_2 \text { and } \alpha_2(x) =
f_2(z,0)
\end{cases}$$

\noindent Likewise, define  $\tilde{\eta}_2:D \times \{0,1 \} \cup
C_2 \to I$ such that:

$$\tilde{\eta}_2(x,t) = \begin{cases}
t & \text{ if } t=0,1 \\
\tau_2(z,t) & \text{ if } x \in Q_2 \text { and } \alpha_2(x) =
f_2(z,0)
\end{cases}$$

\noindent Since $X$ is an ANR, $\tilde{g}_2$ extends to a map $g_2:D
\times I \to X$ and $\tilde{\eta}_2$ extends to a map $\eta_2:D
\times I \to I$. Note that we have used sufficient care in our
construction so that:

\begin{enumerate}
\item  $({g}_1 |_{B_1}, {\eta}_1 |_{B_1})$ and $({g}_2 |_{C_2},
{\eta}_2 |_{C_2})$ are disjoint topographical map pairs; and
\item $({g}_1 |_{C_1}, {\eta}_1 |_{C_1})$ and $({g}_2 |_{B_2},
{\eta}_2 |_{B_2})$ are disjoint topographical map pairs.
\end{enumerate}
\noindent By the $\mathcal{D}_c \times \mathcal{D}$ DTP*, $(g_i,\eta_i)$
can be approximated by disjoint topographical map pairs.  Hence we
also have that
\begin{enumerate}
\item[(3)] $({g}_i|_{C_i},{\eta}_i|_{C_i})$ can be approximated by disjoint
topographical map pairs.
\end{enumerate}
\noindent Therefore, by applying the Special DTP Extension Theorem,
there are disjoint topographical map pairs  $(h_i, \mu_i)$ that are
approximations of $({g}_i |_{C_i}, {\eta}_i|_{C_i})$ so that $(h_i
|_{B_i}, \mu_i |_{B_i})=({g}_i |_{B_i}, {\eta}_i |_{B_i})$. This
determines approximations $(f'_i, \tau'_i)$ of $(f_i, \tau_i)$ that
are disjoint topographies. In particular, $(f'_i(z,t),\tau'_i(z,t)) =
(h_i(x,t),\mu_i(x,t))$, where $f_i(z,0)= \alpha_i(x)$.

\end{proof}

\begin{remark}
Note that it is an easy matter to show that conditions (a)-(d) imply
the DAP. Given two singular arcs in $X$, use these paths to define
constant path homotopies and apply any one of the conditions to
approximate by disjoint topographical map pairs in the case of
(a)-(c) or disjoint concordances in the case of (d). The end levels
provide the disjoint approximations. The equivalence of (e) with
(a)-(d) does require the DAP.
\end{remark}

\section{Crinkled Ribbons Properties} \label{Crinkled}

Recall that, given $k\geq 0$, a subset $Z\subset X$ of space~$X$ is
said to be {\it locally $k$-coconnected}  ($k$-LCC) if for every
point $x\in X$ and every neighbourhood $U\subset X$ of~$x$, there
exists a neighbourhood $V\subset U$ of~$x$ such that the
inclusion-induced homomorphism $\pi_k(V\setminus
Z)\to\pi_k(U\setminus Z)$ is trivial. Also recall the following
useful Proposition (see \cite[Corollary 26.2A]{Daverman book}):

\begin{prop}
Each $k$-dimensional closed subset $A$ of a generalized $n$-manifold
$X$, where $k \leq n-2$, is $0$-LCC.
\end{prop}

We are now ready to define the ribbons properties.

\begin{defn}
A generalized $n$-manifold $X$ has the \emph{crinkled ribbons
property (CRP)} provided that any constant homotopy $f: K \times I
\to X$, where $K$ is a $1$-complex can be approximated by a map
$f':K \times I \to X$ so that:
\begin{enumerate}
\item $f'(K \times \{0\}) \cap f'(K \times \{1\}) = \emptyset$; and
\item dim$(f'(K \times I))\leq n-2$.
\end{enumerate}
\end{defn}

\begin{thm} \label{crp}
If $X$ is a resolvable generalized $n$-manifold, $n \geq 4$, with
the crinkled ribbons property, then $X$ has the $\mathcal{K}_c \times
\mathcal{K}$ DTP*.
\end{thm}

\begin{proof}  Let  $(f_1, \tau_1) \in \mathcal{K}_c$ and $(f_2, \tau_2)
\in \mathcal{K}$.  Apply the hypothesis of the theorem to find
$f'_1: K_1 \times I \to X$ so that $f_1'(K_1 \times \{0\} ) \cap
f_1'(K_1 \times \{1 \} ) = \emptyset$, and dim$(f_1'(K_1 \times I))
\leq n-2$.  It follows that $f_1'(K_1 \times I)$ is $0$-LCC in $X$.
Let $A_0 = f_1'(K_1 \times \{0\} )$ and $A_1 = f_1'(K_1 \times \{1\}
)$. Define $\tau'_1: K_1 \times I \to I$ so that:
$$ \tau'_1 (x,t) =
\dfrac{d(f'_1(x,t),A_0)}{d(f'_1(x,t),A_0)+d(f'_1(x,t),A_1)}.$$ Apply
the $0$-LCC condition to approximate $f_2$ by $f'_2:K_2 \times I \to
X$ so that $$f'_2(K_2 \times [\mathbb{Q}\cap I] \cup {\mathbb{K}_2}
\times I) \cap f'_1(K_1 \times I) = \emptyset,$$ where
${\mathbb{K}_2}$ is a countable dense set in $K_2$ containing the
vertex set.  Then $(f'_2)^{-1}(f'_1(K_1 \times I ))$ is closed
{0}-dimensional set $Z$. Approximate $\tau_2$ so that $\tau'_2$ is
1-1 on $Z$.  Then PPIN$((f'_1, \tau'_1), (f'_2, \tau'_2))$ is a
closed $0$-dimensional set.  By the Reparametrization Lemma, $(f'_i,
\tau'_i)$ can be approximated by disjoint topographical map pairs.
\end{proof}

\begin{cor}
If $X$ is a resolvable generalized $n$-manifold, $n \geq 4$, with
the crinkled ribbons property, then $X \times \mathbb{R}$ has the
disjoint disks property.
\end{cor}

\begin{proof}
Follows directly from Theorem \ref{crp} and the Equivalence Theorem.
\end{proof}

\begin{defn}
A generalized $n$-manifold $X$ has the \emph{twisted crinkled
ribbons property (CRP-T)} provided that any constant homotopy $f: D
\times I$ can be approximated by a map $f':D \times I$ so that:
\begin{enumerate}
\item $f'(D \times \{0\} ) \cap f'(D \times \{1 \} )$ is a finite set of
points; and
\item dim$(f'(D \times I))\leq n-2$.
\end{enumerate}
\end{defn}

\begin{thm} \label{crpt}
If $X$ is a generalized $n$-manifold of dimension $n \geq 4$ having
the twisted crinkled ribbons property and the property that points
are $1$-LCC embedded in $X$, then $X \times \mathbb{R}$ has the
$\mathcal{D}_c \times \mathcal{D}$ DTP*.
\end{thm}

\begin{proof}  It suffices to show that maps in $\mathcal{D}_c
\times \mathcal{D}$ can be approximated by disjoint topographical
map pairs.

Let $(f_1, \tau_1) \in \mathcal{D}_c$ and $(f_2, \tau_2) \in
\mathcal{D}$. Since any generalized manifold of dimension $\geq 3$
has the DAP, we may assume without loss of generality that $(f_i,
\tau_i)$ are disjoint on the end levels, (i.e., $f_1(D \times \{e\})
\cap f_2(D \times \{e\}) = \emptyset$ for $e=0,1$),  and that any
adjustments hereafter are sufficiently small to maintain this
condition.  Apply the hypothesis of the theorem to find $f'_1: D
\times I \to X$ so that $f_1'(D \times \{0\} ) \cap f_1'(D \times
\{1 \} )$ is a finite set of points $P$, and dim$(f_1'(D \times I))
\leq n-2$. It follows that $f_1'(D \times I)$ is $0$-LCC in $X$. We
may also apply the hypothesis that points are $1$-LCC embedded in
$X$ and assume without loss of generality that $f_2(D \times I) \cap
P = \emptyset$.

Choose $\zeta > 0$ so that $d(f_2(D \times I), P ) >
\zeta$.  Let $A_0 = f_1'(D \times \{0\} )- N(P,\zeta)$ and $A_1 =
f_1'(D \times \{1\} )- N(P,\zeta)$. Define
$$\tau'_1:  D \times \{0,1
\} \cup  (D \times I - (f'_1)^{-1} (N(P,\zeta)))  \to I$$ so that:
$$ \tau'_1 (x,t) =\left\{
                    \begin{array}{ll}
                      e & \hbox{ if } t=e \\
                      \dfrac{d(f'_1(x,t),A_0)}{d(f'_1(x,t),A_0)+d(f'_1(x,t),A_1)} & \hbox{ otherwise.}
                    \end{array}
                  \right.$$
Since $D \times I $ is an AR, $\tau'_1$ may be extended to all of $D
\times I$.  Apply the $0$-LCC condition to approximate $f_2$ by
$f'_2:D \times I \to X$ so that:
$$f'_2(D \times [\mathbb{Q}\cap I] \cup [\mathbb{Q}\cap D] \times I)
\cap f'_1(D \times I) = \emptyset.$$ Then $(f'_2)^{-1}(f'_1(D \times
I ))$ is a
closed {0}-dimensional set $Z$. Approximate $\tau_2$ so
that $\tau'_2$ is 1-1 on $Z$.  Then PPIN$((f'_1, \tau'_1), (f'_2,
\tau'_2))$ is a closed $0$-dimensional set.  By the
Reparametrization Lemma, $(f'_i, \tau'_i)$ can be approximated by
disjoint topographical map pairs.
\end{proof}

\begin{cor}
If $X$ is a generalized $n$-manifold of dimension $n \geq 4$ having
 the twisted crinkled ribbons property and the property that points
are $1$-LCC embedded in $X$, then $X \times \mathbb{R}$ has the
disjoint disks property.
\end{cor}

\begin{proof}
The assertion follows directly from Theorem \ref{crpt} and the Equivalence
Theorem.  Note that the condition that points are $1$-LCC embedded
implies the $(0,2)$-DDP.
\end{proof}

\begin{remark}
Not all generalized manifolds of dimension $n \geq 4$ have the
property that points are $1$-LCC embedded.  For example, the
Daverman-Walsh $2$-ghastly spaces are resolvable generalized
manifolds that do not have the $(0,2)$-DDP, and hence cannot satisfy
the condition that points are $1$-LCC embedded
\cite{Daverman-Walsh}.
\end{remark}

The following corollary was also proved in \cite{Daverman book,
Daverman 2} by using shrinking techniques.  This is the first time
general position arguments have been applied to this setting.

\begin{cor}
If $X$ is a resolvable generalized locally spherical $n$-manifold,
$n\geq 4$, then $X$ is a codimension one manifold factor.
\end{cor}

\begin{proof}
The locally spherical condition implies the twisted crinkled ribbons
property.  To see this, let $f:D \times I$ be a constant homotopy.
Cover the image of $f$ by small neighborhoods $B_1, B_2, \ldots B_n$
so that $\partial B_i$ is an embedded $(n-1)$-sphere.  Approximate
$f$ by a constant path homotopy $f': D \times I \to \bigcup
\partial B_i $.  Without loss of
generality we may assume that
there are $t_i \in D$ such that $0 = t_0 < t_1 < \ldots <
t_{n-1} < t_n = 1$ and  $f'([t_{i-1}, t_i]\times I ) \subset
\partial B_n$. Since $\partial B_i$ is an $(n-1)$-sphere, $f'$ can be approximated by $f'': D \times I \to
\bigcup
\partial B_i$ such that the
restriction of $f''$ to $\bigcup(t_{i-1}, t_i) \times I$ is an
embedding and $f'' = f'$ on $\{t_0, t_1, \ldots  t_n) \times I$.
Then $f''$ is the desired approximation of $f$.
\end{proof}

\section{Fuzzy Ribbons property} \label{Fuzzy}

The fuzzy ribbons property is the most remarkable generalization of
the disjoint homotopies techniques.  In particular the fuzzy ribbons
property is a generalization of the method of $\delta$-fractured
maps. Recall that

\begin{defn}
A map $f:D \times I \to X$ is said to be
\emph{$\delta$-fractured} over a
map $g:D \times I \to X$ if there are pairwise disjoint balls $B_1, B_2,
\ldots, B_m$ in $D \times I$ such that:
    \begin{enumerate}
        \item $diam(B_i) < \delta$;
        \item $f^{-1}(im(g)) \subset \bigcup_{i=1}^m int(B_i)$; and
        \item $diam(g^{-1}(f(B_i))) < \delta$.
    \end{enumerate}
\end{defn}

\noindent However, because of the freedom in defining the level map
to obtain the DTP* conditions, we need no longer require
$\delta$-control.  The analogous definition in the setting of
topographical map pairs is therefore:

\begin{defn}
Let $(f_i, \tau_i) \in \mathcal{K}$ be such that $f_i: K_i \times I \to
X$ and $\tau_i: K_i \times I \to I$.  Then $(f_2, \tau_2)$ is said
to be {\it fractured} over a topographical map pair $(f_1, \tau_1)$ if
there are disjoint balls $B_1, B_2, \ldots, B_m$ in $K_2 \times I$
such that:
    \begin{enumerate}
        \item $f_2^{-1}(im(f_1)) \subset \bigcup_{j=1}^m int(B_i)$; and
        \item $\tau_1 \circ f_1^{-1}\circ f_2(B_i) \ne I$.
    \end{enumerate}
\end{defn}

\noindent We are now ready to define the fuzzy ribbons property:

\begin{defn}
A space $X$ has the \emph{fuzzy ribbons property (FRP)} provided
that for any topographical map pairs, $(f_1, \tau_1) \in
\mathcal{K}_c$ and $(f_2, \tau_2) \in \mathcal{K}$, and $\varepsilon
> 0$ there are maps $\tau'_i$ and $\varepsilon$-approximations $f'_i$ of  $f_i$  so that
$(f'_2, \tau'_2)$ is fractured over $(f'_1, \tau'_1)$.
\end{defn}

\begin{thm} \label{fractured theorem}
If a space $X$ is an ANR with the DAP having the fuzzy ribbons property,
then $X$ has the $\mathcal{K}_c \times \mathcal{K}$ DTP*.
\end{thm}

\begin{proof}
Let $(f_1, \tau_1) \in \mathcal{K}_c$ and $(f_2, \tau_2) \in
\mathcal{K}$ such that $f_i: K_i \to X$.  Using the DAP we may
assume without loss of generality that $f_1(K_1 \times I) \cap
f_2(K_2 \times \{0,1\} \cup K^{(0)}_2 \times I) = \emptyset$. Apply
the fuzzy ribbons property to obtain maps $\tau'_i$ and
approximations $f'_i$ of $f_i$ so that $(f'_2, \tau'_2)$ is
fractured over $(f'_1, \tau'_1)$. The approximations of $f_i$ should
be sufficiently small so that $f'_1(K_1 \times I) \cap f'_2(K_2
\times \{0,1\} \cup K^{(0)}_2 \times I) = \emptyset$. Then there are
disjoint balls $B_1, B_2, \ldots, B_m$ in $K_2 \times I - K_2 \times
\{0,1\} \cup K^{(0)}_2 \times I$ such that:
    \begin{enumerate}
        \item $(f'_2)^{-1}(im(f'_1)) \subset \bigcup_{j=1}^m int(B_i)$; and
        \item $\tau'_1 \circ (f'_1)^{-1}\circ f'_2(B_i) \ne I$.
    \end{enumerate}
For each $j = 1, \ldots, m$, choose $t_j \in I - \tau'_1 \circ
(f'_1)^{-1}\circ f'_2(B_i)$.  Now define $\tau''_2: K_2 \times \{0,1
\} \cup K^{(0)}_2 \times I \cup (\bigcup B_j) \to I$ so that
$$ \tau''_2(x,t) = \left\{
                                                 \begin{array}{ll}
                                                   t & \hbox{ if } t=0,1 \hbox{ or } x \in K^{(0)}_2\\
                                                   t_j & \hbox{ if
} (x,t) \in B_j
                                                 \end{array}
                                               \right.$$
Extend $\tau''_2$ to $K_2 \times I$.  Then $(f'_1, \tau'_1)$ and
$(f'_2, \tau''_2)$ are the desired disjoint topographical map pairs.
\end{proof}

\begin{cor}
If a space $X$ is an ANR with the the FRP, then $X \times \mathbb{R}$
has the DDP.
\end{cor}

\begin{proof}
The DAP follows from the FRP.  The rest follows directly from Theorem
\ref{fractured theorem}
and
the Equivalence Theorem.
\end{proof}

\begin{remark}
Certain $2$-ghastly spaces satisfy the FRP, such as those discussed in \cite{Halverson 2}. The same type of arguments apply, however less attention to control is needed to satisfy the FRP.
\end{remark}

\section{Epilogue}
The DTP* properties presented in this paper are not only more
versatile in detecting codimension one manifold factors, they also
provide a characterization of such spaces.  The ribbons properties
represent practical applications of these properties.  Further interesting questions that may be investigated using the DTP* or ribbons properties include:

\begin{que}
   If $G$ is an $(n-2)$-dimensional cell-like decomposition of an
$n$-manifold $M$, where $n \geq 4$, is $M/G$ a codimension one
manifold factor?
\end{que}

\begin{que}
Is every Busemann $G$-space of dimension $n \geq 5$ a manifold? Equivalently,
are small metric spheres in these spaces codimension one manifold
factors?
\end{que}

\begin{que}
Is every finite-dimensional resolvable generalized manifold of
dimension $n\geq 4$ a codimension one manifold factor?
\end{que}

\section*{Acknowledgments}

This research was supported by the Brigham Young University
Special Topology Year 2008-2009 Fund and
the Slovenian Research Agency grants
BI-US/07-08/029,
P1-0292-0101, and
J1-9643-0101.
We thank the referee for several comments and suggestions.


\begin{thebibliography}{99}

\bibitem{Daverman book}
R. J. Daverman,
Decompositions of Manifolds. Pure and Applied Mathematics, 124.
Academic Press, Inc., Orlando, FL, 1986.

\bibitem{Daverman 2}
R. J. Daverman,
{\it Singular regular neighborhoods and local flatness in codimension one,}
Proc. Amer. Math. Soc. {\bf 57} (1976), 357--362.

\bibitem{Daverman}
R. J. Daverman,
{\it  Embedding phenomena based upon decomposition theory:
locally spherical but wild codimension one spheres,}
Proc. Amer. Math. Soc. {\bf 90}:1 (1984), 139--144.

\bibitem{Daverman-Halverson 2}
R. J. Daverman and D. M. Halverson,
{\it Path concordances as detectors of codimension-one manifold factors},
Exotic Homology Manifolds, Oberwolfach 2003,  Geom. Topol. Monogr. {\bf 9},
Geom. Topol. Publ., Coventry, 2006, pp. 7--15.

\bibitem{Daverman-Halverson}
R. J. Daverman and D. M. Halverson,
{\it The cell-like approximation theorem in dimension $n = 5$,}
Fund. Math. {\bf 197}  (2007), 81--121.

\bibitem{Daverman-Walsh}
R. J. Daverman and J. J. Walsh,
{\it A ghastly generalized $n$-manifold},
Illinois J. Math. {\bf 25}:4 (1981), 555--576.

\bibitem{Edwards 1}
R. D. Edwards,
{\it Demension theory, I},
Geometric Topology, Proc. Conf., Park City, Utah, 1974,
Lecture Notes in Math. {\bf 438}, Springer, Berlin, 1975, pp. 195--211.

\bibitem{Edwards 2}
R. D. Edwards,
{\it The topology of manifolds and cell-like maps},
Proc. Int. Congr. of Math., Helsinki, 1978,
Acad. Sci. Fennica, Helsinki, 1980, pp. 111--127

\bibitem{Halverson 1}
D. M. Halverson,
{\it Detecting codimension one
manifold factors with the disjoint homotopies property,}
Topology Appl. {\bf 117}:3 (2002), 231--258.

\bibitem{Halverson 2}
D. M. Halverson,
{\it $2$-ghastly spaces with
the disjoint homotopies property: the method of fractured maps,}
Topology Appl. {\bf 138}:1-3 (2004), 277--286.

\bibitem{Halverson 3}
D. M. Halverson,
{\it Detecting codimension one manifold factors with 0-stitched disks},
Topology Appl. {\bf 154}:9  (2007), 1993--1998.

\bibitem{HaRe}
D. M. Halverson and D. Repov\v{s},
{\it The Bing-Borsuk and the Busemann Conjectures},
Math. Comm. {\bf 13}:2 (2008), 163--184.

\bibitem{Repovs1}
D. Repov\v{s},
{\it The recognition problem for topological manifolds},
Geometric and Algebraic Topology, J. Krasinkiewicz, S. Spie\.{z}, and H. Toru\,{n}czyk, Eds.,
PWN, Warsaw 1986, pp. 77--108.

\bibitem{Repovs2}
D. Repov\v{s},
{\it  Detection of higher dimensional topological manifolds among topological spaces},
Giornate di Topologia e Geometria Delle Variet\`{a}, Bologna 1990, M. Ferri, Ed.,
Univ. degli Studi di Bologna 1992, pp. 113--143.

\bibitem{Repovs3}
D. Repov\v{s},
{\it The recognition problem for topological manifolds: A survey,}
Kodai Math. J. {\bf 17}:3 (1994), 538--548.

\end{thebibliography}
\end{document}